\documentclass[11pt, a4paper, leqno]{amsart}
\usepackage{amssymb}
\author{Raf Cluckers$^*$)}
\thanks{$^*$) Research Assistant of the Fund for Scientific Research --
 Flanders (Belgium)(F.W.O.)}
\date{}
\title{Classification of semi-algebraic $p$-adic sets up to semi-algebraic bijection.}

\newtheorem{theorem}{Theorem}[section]
\newtheorem{lemma}[theorem]{Lemma}
\newtheorem{proposition}[theorem]{Proposition}
\newtheorem{cor}[theorem]{Corollary}

\theoremstyle{definition}
\newtheorem{definition}[theorem]{Definition}

\theoremstyle{remark}
\newtheorem{remark}[theorem]{Remark}

\newcommand{\Q}{\mathbb{Q}}
\newcommand{\F}{\mathbb{F}}
\newcommand{\Z}{\mathbb{Z}}
\newcommand{\N}{\mathbb{N}}
\newcommand{\R}{\mathbb{R}}

\DeclareMathOperator{\sq}{\square}
\begin{document}
\maketitle
\begin{abstract}
We prove that two infinite $p$-adic semi-algebraic sets are isomorphic
(i.e. there exists a semi-algebraic bijection between them) if and only if
they have the same dimension.
\end{abstract}

In real semi-algebraic geometry (as opposed to $p$-adic
semi-algebraic geometry) the following classification is
well-known \cite{vdD}:
 \begin{quote}\textit{There exists a real semi-algebraic bijection
between two real semi-algebraic sets if and only if they have the
same dimension and Euler characteristic.}
\end{quote}

More generally L. van den Dries \cite{vdD} gave such a
classification for o-minimal expansions of the real field, using
the dimension and Euler characteristic as defined for $o$-minimal
structures. Since the semi-algebraic Euler characteristic
 $\chi$ is in fact the canonical map from
the real semi-algebraic sets onto the Grothendieck ring (see
\cite{CH}) of $\R$ (which is $\Z$), we see that the isomorphism
class of a real semi-algebraic set only depends on its image in
the Grothendieck ring and its dimension.
\par
In this paper we treat the $p$-adic analogue of this
classification. The
 Grothendieck ring of $\Q_p$ is recently proved to be trivial by D.
Haskell and the author \cite{CH}, so the analogue of the real case
is a classification of the $p$-adic semi-algebraic sets up to
semi-algebraic bijection using only the dimension. We give such a
classification for the $p$-adic semi-algebraic sets and for finite
field extensions of $\Q_p$, using explicit isomorphisms of
\cite{CH} and the $p$-adic Cell Decomposition Theorem of J. Denef
\cite{D}. The most difficult part in giving this classification is
to prove that for any semi-algebraic set $X$ there is a finite
partition into semi-algebraic sets, such that each part is
isomorphic to a Cartesian product of one dimensional sets, in
other words semi-algebraic sets have a rectilinearization. Since
all arguments hold also for finite field extensions of $\Q_p$, we
work in this more general setting.
\par
Let $p$ denote a fixed prime number, $\Q_p$ the field of $p$-adic
numbers and $K$ a fixed finite field extension of $\Q_p$. For
$x\in K$ let $v(x)\in\Z\cup\{+\infty\}$ denote the valuation of
$x$. Let $R=\{x\in K\mid v(x)\geq0\}$ be the valuation ring,
$K^\times=K\setminus\{0\}$ and for $n\in\N_0$ let $P_n$ be the set
$\{x\in K^\times\mid\exists y\in K\ y^n=x\}$. We call a subset of
$K^n$ {\em semi-algebraic } if it is a Boolean combination (i.e.
obtained by taking finite unions, complements and intersections)
of sets of the form $\{x\in K^m|f(x)\in P_n\}$, with $f(x)\in
K[X_1,\ldots,X_m]$. The collection of semi-algebraic sets is
closed under taking projections $K^m\to K^{m-1}$, even more: it
consists precisely of Boolean combinations of projections of
affine $p$-adic varieties. Further we have that sets of the form
$\{x\in K^m|v(f(x))\leq v(g(x))\}$ with $f(x),g(x)\in
K[X_1,\ldots,X_m]$ are semi-algebraic (see \cite{D} and \cite{M}).
A function $f:A\to B$ is semi-algebraic if its graph is a
semi-algebraic set; if further $f$ is a bijection, we call $f$
 an {\em isomorphism } and we write $A\cong B$.
\par
 Let $\pi$ be a fixed element of $R$ with
$v(\pi)=1$, thus $\pi$ is a uniformizing parameter for $R$. For a
semi-algebraic set $X\subset K$ and $k>0$ we write
\[
X^{(k)}=\{x\in X|x\not=0 \mbox{ and } v(\pi^{-v(x)}x-1)\geq k,\
x\not=0\},
\]
which is semi-algebraic (see \cite{D}, Lemma 2.1); $X^{(k)}$ consists of
those points $x\in X$ which have a $p$-adic expansion $x=\sum_{i=s}^\infty
a_i\pi^i$ with $a_s=1$ and $a_i=0$ for $i=s+1,\ldots,s+k-1$. By a
\emph{finite partition} of a semi-algebraic set we mean a partition into
finitely many semi-algebraic sets. Let $X\subset K^n$, $Y\subset K^m$ be
semi-algebraic. Choose disjoint semi-algebraic sets $X',Y'\subset K^k$ for
some $k$, such that $X\cong X'$ and $Y\cong Y'$, then we define the
disjoint union of $X$ and $Y$ up to isomorphism as $X'\cup Y'$. In the
introduction of \cite{CH} it is shown that we can take $k=\max(m,n)$, i.e.
we can realize the disjoint union without going into higher dimensional
affine spaces.
\par
We recall some well-known facts.
\begin{lemma}[Hensel]\label{Hensel}
Let $f(t)$ be a polynomial over $R$ in one variable $t$, and let
$\alpha\in R,\ e\in\N$. Suppose that
$f(\alpha)\equiv0\bmod\pi^{2e+1}$ and $
  v(f'(\alpha))\leq e$, where $f'$ denotes the derivative of $f$. Then there
exists a unique $\bar\alpha\in R$ such that $f(\bar\alpha)=0$ and
$\bar\alpha\equiv\alpha\bmod\pi^{e+1}$.
\end{lemma}
\begin{cor}\label{corhensel}
Let $n>1$ be a natural number. For each $k> v(n)$, and $k'=k+ v(n)$ the
function $$ K^{(k)}  \to P_n^{(k')}:x\mapsto x^n $$ is an isomorphism.
\end{cor}
\par
The next theorem gives some concrete isomorphisms between one
dimensional sets.
\begin{proposition}[\cite{CH}, Prop. 2]\label{gring}
\item{(i)} The union of two disjoint
copies of $R\setminus\{0\}$ is isomorphic to $R\setminus\{0\}$.
\item{(ii)}  For each $k>0$ the union of two disjoint
copies of $R^{(k)}$ is isomorphic to $R^{(k)}$.
\item{(iii)} $R\cong R\setminus\{0\}$.
\end{proposition}
We deduce an easy corollary, also consisting of concrete isomorphisms.
\begin{cor}\label{corgring} For each $k$ we have isomorphisms
\item{(i)} $R^{(k)}\cong R\setminus\{0\}$,
\item{(ii)} $R\setminus \{0\}\cong K$.
\end{cor}
\begin{proof} (i) There is a finite partition
$R\setminus\{0\}=\bigcup_\alpha \alpha R^{(k)}$ with $v(\alpha)=0$, say
with $s$ parts. Then $R\setminus\{0\}$ is a fortiori isomorphic to the
union of $s$ disjoint copies of $R^{(k)}$,
which is by Proposition \ref{gring}(ii) isomorphic to $R^{(k)}$.\\
 (ii) The map
\[
 (\{0\}\times R)\cup (\{1\}\times R\setminus\{0\})\to K:\left\{\begin{array}{rcl}
                    (0,x) & \mapsto & x,\\
                    (1,x) & \mapsto & 1/(\pi x),\end{array}\right.\
\]
is a well-defined isomorphism. It follows that $K$ is isomorphic to the
disjoint union of $R$ and $R\setminus\{0\}$. Now use (i) and (iii) of
Proposition~\ref{gring}.
\end{proof}
\par Give $K^m$ the topology induced by the norm
$|x|=\max(|x_i|_p)$ with $|x_i|_p=p^{- v(x_i)}$ for $x=(x_1,\ldots
x_m)\in K^m$. P. Scowcroft and L. van den Dries \cite{SvdD} proved
there exists no isomorphism from an open set $A\subset K^m$ onto
an open set $B\subset K^n$ with $n\not= m$, so we can define the
dimension of semi-algebraic sets as follows.
\begin{definition}[\cite{SvdD}]\textup{
 The \emph{dimension} of a semi-algebraic set $X\not=\phi$ is the greatest
natural number $n$ such that we have a nonempty semi-algebraic
subset $A\subset X$ and an isomorphism from $A$ to a nonempty
semi-algebraic open subset of $K^n$. We put $\dim(\phi)=-1$.
 }\end{definition}
\par
P. Scowcroft and L. van den Dries \cite{SvdD} proved many good
properties of this dimension, for example that it is invariant
under isomorphisms.
\begin{proposition}[\cite{SvdD}]
Let $A$ and $B$ be semi-algebraic sets, then the following is true:
\item{(i)} If $A\cong B$ then $\dim(A)=\dim(B)$,
\item{(ii)} $\dim(A\cup
B)=\max(\dim(A),\dim(B))$.
\item{(iii)} $\dim (A)=0$ if and only if $A$ is finite and
nonempty.
\end{proposition}
We will prove the converse of (i) for infinite semi-algebraic sets.
\begin{lemma}\label{embed}
For any semi-algebraic set $X$ of dimension $m\in\N_0$ there exists a
semi-algebraic injection $X\to K^m$.
\end{lemma}
\begin{proof} By \cite{SvdD},~Cor. 3.1 there is a finite partition of
$X$ such that each part $A$ is isomorphic to a semi-algebraic open
$A'\subset K^k$ for some $k\leq m$. Now realize the disjoint union of the
sets $A'$ without going into higher embedding dimension (see the
introduction).
\end{proof}
\par
We formulate the $p$-adic Cell Decomposition Theorem by J. Denef
\cite{D, D2}, which is the analogue of the real semi-algebraic
Cell Decomposition Theorem.
\begin{theorem}[Cell Decomposition \cite{D,D2}] Let $x=(x_1,\ldots,x_m)$ and
$\hat x=(x_1,\ldots,x_{m-1}),\ m>0$. Let $f_i(\hat x,x_m)$,
$i=1,\ldots,r$, be polynomials in $x_m$ with coefficients which are
semi-algebraic functions from $K^{m-1}$ to $K$. Let $n\in\N_0$ be fixed.
Then there exists a finite partition of $K^m$ into sets $A$ of the form
\[A=\{x\in K^m|\hat x\in D\mbox{ and }  v(a_1(\hat x))\sq_1
 v(x_m-c(\hat x))\sq_2 v(a_2(\hat x))\},\]
 such that
\[f_i(x)=u_i(x)^nh_i(\hat x)(x_m-c(\hat x))^{\nu_i},\mbox{ for each }\ x\in A,
\ i=1,\ldots,r,\] with $u_i(x)$ a unit in $R$ for each $x$,
$D\subset K^{m-1}$ semi-algebraic, $\nu_i\in\N$, $h_i,a_1,a_2,c$
semi-algebraic functions from $K^{m-1}$ to $K$ and $\sq_1$,
$\sq_2$ either $\leq,<,$ or no condition.
\end{theorem}
The next lemma is also due to J. Denef \cite{D2}.
\begin{lemma}[\cite{D2}, Cor. 6.5]\label{order}
Let $b:K^m\to K$ be a semi-algebraic function. Then there exists a finite
partition of $K^m$ such that for each part $A$ we have $e>0$ and
polynomials $f_1, f_2\in R[X_1,\ldots X_m]$ such that
\[
v(b(x))=\frac{1}{e} v(\frac{f_1(x)}{f_2(x)}), \mbox{ for each }\ x\in A,
\]
with $f_2(x)\not=0$ for each $x\in A$.
\end{lemma}
We give an application of the Cell Decomposition Theorem and Lemma
\ref{order}, inspired by similar applications in \cite{D2}. For details of
the proof we refer to the proof of \cite{D2}, Thm. 7.4. By $\lambda P_n$
with $\lambda=0$ we mean $\{0\}$.
\begin{lemma}\label{partition}
Let $X\subset K^m$ be semi-algebraic and $b_j:K^m\to K$ semi-algebraic
functions for $j=1,\ldots,r$. Then there exists a finite partition of $X$
s.t. each part $A$ has the form
\begin{eqnarray*}
A & = & \{x\in K^m|\hat x\in D, \  v(a_1(\hat
    x))\sq_1 v(x_m-c(\hat x))\sq_2 v(a_2(\hat x)),\\
 & & \qquad\qquad\qquad x_m-c(\hat x)\in \lambda P_n\},\end{eqnarray*}
and such that for each $x\in A$ we have
\[ v(b_j(x))=\frac{1}{e_j} v((x_m-c(\hat x))^{\mu_j}d_j(\hat x)),\]
with $\hat x=(x_1,\ldots, x_{m-1})$, $D\subset K^{m-1}$
semi-algebraic, $e_j>0$, $\mu_j\in\Z$, $\lambda\in K$, $c,a_i,d_j$
semi-algebraic functions from $K^{m-1}$ to $K$ and $\sq_i$ either
$<,\leq$ or no condition.
\end{lemma}
\begin{proof}
By Lemma \ref{order} we have a finite partition of $X$ such that for each
part $A_0$ we have $e_j>0$ and polynomials $g_j, g'_j\in
R[X_1,\ldots,X_m]$ with
\[
v(b_j(x))=\frac{1}{e_j} v(\frac{g_j(x)}{g'_j(x)}), \mbox{ for each }x\in
A_0,\ j=1,\ldots,r.
\]
Let $f_i$ be the polynomials which appear in a description of
$A_0$ as a Boolean combination of sets of the form $\{x\in K^m\mid
f(x)\in P_n\}$.. Apply now the Cell Decomposition Theorem as in
the proof of \cite{D2}, Thm. 7.4 to the polynomials $f_i,g_j$ and
$g_j'$ to obtain the lemma.
\end{proof}
\par
The proof of the next proposition is an application of both the
Cell Decomposition Theorem and some hidden Presburger arithmetic
in the value group of $K$; it is the technical heart of this
paper. If $l=0$ then $\prod_{i=1}^l R^{(k)}$ denotes the set
$\{0\}$.
\begin{definition}\textup{
We say that a semi-algebraic function $f:B\to K$ satisfies condition
(\ref{condition}) (with constants $e,\mu_i,\beta$) if we have constants
$e\in\N_0,\mu_i\in\Z,\beta\in K$ such that each $x=(x_i)\in B$ satisfies
\begin{equation}\label{condition}
 v(f(x))=\frac{1}{e} v(\beta\prod_i x_i^{\mu_i}).
 \end{equation}
 }\end{definition}
\begin{proposition}[Rectilinearization]\label{resolution}
 Let $X$ be a semi-algebraic set and $b_j:X\to K$ semi-algebraic functions for
$j=1,\ldots,r$. Then there exists a finite partition of $X$ such that for
each part $A$ we have constants $l\in\N$, $k\in\N_0$, $\mu_{ij}\in\Z$,
$\beta_j\in K$, and an isomorphism
\[f:\prod_{i=1}^l R^{(k)}\to A,\]
such that for each $x=(x_1,\ldots,x_l)\in\prod_{i=1}^l R^{(k)}$ we have
\[ v(b_j\circ f(x))= v(\beta_j\prod_{i=1}^{l}x_i^{\mu_{ij}}).\]
\end{proposition}
\begin{proof}
We work by induction on $m=\dim(X)$. Let $\dim(X)=1$ and $b_j:X\to
K$ semi-algebraic functions, $j=1,\ldots,r$. By Lemma \ref{embed}
we may suppose that $X\subset K$. We reduce first to the case that
$X$ and $b_j$ have the special form (\ref{eq dim1}) (see below).
By Lemma \ref{partition} there is a partition such that each part
$A$ is either a point or of the form
\[A=\{x\in
K| v(a_1)\sq_1 v(x-c)\sq_2 v(a_2),\
 x-c\in \lambda P_n\},\]
and such that for each $x\in A$ we have $ v(b_j(x))=\frac{1}{e_j}
v(\beta_j(x-c)^{\mu_j}),$ with $a_i,c,\lambda,\beta_j\in K$,
$e_j>0$ and $\mu_j\in\Z$. We may assume that $\lambda\not=0$,
$a_1\not=0\not=a_2$, $\sq_i$ is either $\leq$ or no condition and
since the translation
\[
\{x\in K| v(a_1)\sq_1 v(x)\sq_2 v(a_2),\ x\in \lambda P_n\}\to
A:x\mapsto x+c
\]
is an isomorphism, we may also assume that $c=0$. If both $\sq_1$
and $\sq_2$ are no condition we can partition $A$ into parts
$\{x\in A|0\leq v(x)\}$ and $\{x\in A| v(x)\leq-1\}$. It follows
that if $\sq_1$ is no condition we may suppose that $\sq_2$ is
$\leq$, then we can apply the isomorphism
 \[
 \{x\in K| v(\frac{1}{a_2})\leq v(x),\ x\in
\frac{1}{\lambda}P_n\}\to A:
 x\mapsto\frac{1}{x},
 \]
 and replace $\mu_j$ by $-\mu_j$.
This shows we can reduce to the case that $X$ has the form
\begin{equation}\label{eq dim1}
X=\{x\in K| v(a_1)\leq v(x)\sq_2 v(a_2),\ x\in\lambda P_n\},
\end{equation}
with $a_1\not=0\not=a_2$, $\lambda\not=0$, $\sq_2$ either $\leq$
or no condition and $v(b_j(x))=\frac{1}{e_j} v(\beta_j
x^{\mu_j})$ for each $x\in X$.\\

 \textbf{Case 1:} $\sq_2$ is $\leq$ (in equation (\ref{eq dim1})).\\
By Hensel's Lemma we can partition $X$ into finitely many parts of
the form $y+\pi^sR$ for some fixed $s>v(a_2)$ and with $v(a_1)\leq
v(y)\leq v(a_2)$ for each $y$. For each such part there is a
finite partition
$y+\pi^sR=\bigcup_{\gamma\in\Gamma}A_\gamma\cup\{y\}$, with
$A_\gamma=y+\pi^s\gamma R^{(1)}$ and $v(\gamma)=0$ for each
$\gamma$. The functions $f_\gamma:R^{(1)}\to A_\gamma:x\mapsto
y+\pi^s\gamma x$ are isomorphisms which satisfy $v(b_j\circ
f_\gamma(x))=\frac{1}{e_j}v(\beta_jy^{\mu_j})$ for all $x\in
R^{(1)}$. This last expression is independent of $x$, so there
exists $\beta_j'\in K$ such that $v(b_j\circ f_\gamma(x)) =
v(\beta'_j)$ for all $x\in R^{(1)}$. This proves Case 1.\\

 \textbf{Case 2:} $\sq_2$ is no condition (in equation (\ref{eq dim1})).\\
The map
\[f_1:R\cap\lambda' P_n\to X:x\mapsto a_1x,\]
with $\lambda'=\lambda/a_1$ is an isomorphism. Let $n'$ be a
common multiple of $e_1,\ldots,e_r$ and $n$. Choose $k> v(n')$ and
put $k'=k+ v(n')$. Let  $R\cap\lambda' P_n=\bigcup_\gamma
B_\gamma$ be a finite partition, with $B_\gamma=\gamma(R\cap
P_{n'}^{(k')})$ and $0\leq v(\gamma)<n'$. Now we have that the map
$f_\gamma:R^{(k)}\to B_\gamma:x\mapsto\gamma x^{n'}$ is an
isomorphism by corollary \ref{corhensel}. Let $g_\gamma$ be the
semi-algebraic function $f_1\circ f_\gamma$, which is an
isomorphism from $R^{(k)}$ onto a semi-algebraic set
$A_\gamma\subset X$. The sets $A_\gamma$ form a finite partition
of $X$. Put $\mu'_j=\mu_jn'/e_j$, then we have for each $x\in
R^{(k)}$ that
\[
v(b_j\circ g_\gamma(x))=\frac{1}{e_j} v(\beta_j(a_1\gamma
x^{n'})^{\mu_j})=\frac{1}{e_j}v(\beta_j(a_1\gamma)^{\mu_j})+v(x^{\mu_j'})
\]
\[
= v(\beta'_j)+v(x^{\mu'_j})=v(\beta'_jx^{\mu'_j}),
\]
with $\beta'_j$ in $K$ such that $ v(\beta'_j)=\frac{1}{e_j}
v(\beta_j(a_1\gamma)^{\mu_j})$. This proves case 2.
\par
Now let $\dim(X)=m>1$ and let $b_j:X\to K$ be semi-algebraic
functions, $j=1,\ldots,r$. By Lemma \ref{embed} we may suppose
that
$X\subset K^m$.\\

 \textbf{Claim.}\textit{ We can
partition $X$ such that for each part $A$ we have an isomorphism of the
form $f:D_1\times D_{m-1}\to A$, with $D_1\subset K$ and $D_{m-1}\subset
K^{m-1}$ semi-algebraic, such that the functions $b_j\circ f$ satisfy
condition (\ref{condition}), i.e. there are constants
$e_j\in\N_0,\mu_{ij}\in\Z,\beta_j\in K$ such that each $x=(x_i)\in
D_1\times D_{m-1}$ satisfies
\[
 v(b_j\circ f(x))=\frac{1}{e_j} v(\beta_j\prod_i x_i^{\mu_{ij}}).
 \]}

If the claim is true, we can apply the induction hypotheses once to $D_1$
and the functions $x_1\mapsto x_1^{\mu_{1j}}$ and once to $D_{m-1}$ and
the functions $(x_2,\ldots,x_m)\mapsto\beta_j\prod_{i=2}^mx_i^{\mu_{ij}}$
for $j=1,\ldots, r$. It follows easily that we can partition $X$ such that
for each part $A$ there is an isomorphism $f:\prod_iR^{(k)}\to A$ such
that all $f\circ b_j$ satisfy condition (\ref{condition}) with constants
$e'_j,\mu'_{ij}$ and $\beta'_j$. Now we can proceed as in Case 2 for $m=1$
to make all $e'_j\in\N_0$ occurring
in condition (\ref{condition}) equal to 1. The proposition follows now immediately.\\

\textbf{Proof of the claim.} First we show we can reduce to the
case described in equation (\ref{eq dim>1}) below. Using Lemma
\ref{partition} and its notation, we find a finite partition of
$X$ such that each part $A$ has the form
\begin{eqnarray*}
A & = & \{x\in K^m|\hat x\in D,\  v(a_1(\hat
    x))\sq_1 v(x_m-c(\hat x))\sq_2 v(a_2(\hat x)),\\
 & &\qquad\qquad\qquad x_m-c(\hat x)\in \lambda P_n\},
\end{eqnarray*}
 and such that for each $x\in A$ we have $ v(b_j(x))=\frac{1}{e_j}
v((x_m-c(\hat x))^{\mu_{mj}}d_j(\hat x))$, with $\mu_{mj}\in\Z$. Similar
as for $m=1$, we may suppose that $c(\hat x)=0$ for all $\hat x$. Apply
now the induction hypotheses to the set $D\subset K^{m-1}$ and the
functions $a_1,a_2,d_j$. We find a finite partition of $A$ such that for
each part $A'$ we have an isomorphism $f:B\to A'$, where $B$ is a set of
the form
\[
 B=\{x\in K^m|\hat x\in D',\  v(\alpha_1\prod_{i=1}^lx_i^{\eta_i})
\sq_1 v(x_m)\sq_2 v(\alpha_2\prod_{i=1}^lx_i^{\varepsilon_i}),\
x_m\in\lambda P_n\},
\]
with $D'=\prod_{i=1}^lR^{(k)}$, $l\leq m-1$, such that each
$b_j\circ f$ satisfies condition (\ref{condition}). We will
alternately partition further and apply isomorphisms to the parts
which compositions with $b_j$ will always satisfy condition
(\ref{condition}). By the induction hypotheses we may suppose that
$\lambda\not=0$ and $\dim(D')=m-1$, i.e.
$D'=\prod_{i=1}^{m-1}R^{(k)}$. Analogously as for $m=1$ we may
suppose that $\alpha_1\not=0\not=\alpha_2$, $\sq_2$ is either
$\leq$ or no condition and $\sq_1$ is the symbol $\leq$ (possibly
after partitioning or applying
$x\mapsto(x_1,\ldots,x_{m-1},1/x_m)$).
\par
Choose $\bar k>v(n)$ and put $k'=\bar k+v(n)$. We may suppose that $k'>k$,
so we have a finite partition $B=\bigcup_\gamma B_\gamma$ with
$\gamma=(\gamma_1,\ldots,\gamma_m)\in K^m$, $0\leq v(\gamma_i)<n$ and
$B_\gamma=\{x\in B|x_i\in\gamma_i P_n^{(k')}\}$. Now we have isomorphisms
\[
f_\gamma:C_\gamma\to
B_\gamma:x\mapsto(\gamma_1x_1^n,\ldots,\gamma_mx_m^n),
\]
with
\[
C_\gamma=\{x\in (\prod_{i=1}^{m-1}R^{(\bar k)})\times K^{(\bar
k)}|v(\alpha'_1\prod_{i=1}^{m-1}x_i^{\eta_i}) \leq v(x_m)\sq_2
v(\alpha'_2\prod_{i=1}^{m-1}x_i^{\varepsilon_i})\},
\]
for appropriate choice of $\alpha_i'\in K$. Put
$\nu_i=\varepsilon_i-\eta_i$, $\beta=\alpha'_2/\alpha'_1$, then we have
the isomorphism
\[
\begin{array}{c}
\{x\in \prod_{i=1}^mR^{(\bar k)}|v(x_m)\sq_2
 v(\beta\prod_{i=1}^{m-1} x_i^{\nu_i})\} \to C_\gamma:\\
 x\mapsto(x_1,\ldots,x_{m-1},\alpha'_1x_m\prod_{i=1}^{m-1}x_i^{\eta_i}).
\end{array}
\]
\par
If $\sq_2$ is no condition, the claim is trivial. It follows that
we can reduce to the case that we have an isomorphism
\begin{equation}\label{eq dim>1}
 f:E=\{x\in \prod_{i=1}^mR^{(\bar k)}|v(x_m)\leq
 v(\beta\prod_{i=1}^{m-1} x_i^{\nu_i})\}\to X
\end{equation}
with $\beta\not=0$, $\bar k>0$, and $\nu_i\in\Z$, such that each
$b_j\circ f$  satisfies condition (\ref{condition}).
\par
Suppose we are in the case described in (\ref{eq dim>1}). If
$\nu_i\leq0$ for $i=1,\ldots,m-1$ then we have a finite partition
$E=\cup_sE^{=s}$, with $s\in\{0,1,\ldots,v(\beta)\}$ and
$E^{=s}=\{x\in E| v(x_m)=s\}$. Also,
$E^{=s}=\{(x_1,\ldots,x_{m-1})|\exists x_m\  (x_1,\ldots,x_m)\in
E^{=s}\}\times \{x_m\in R^{(\bar k)}|v(x_m)=s\}$ and the claim
follows.
\par
Suppose now that $\nu_1>0$ in (\ref{eq dim>1}). First we prove the
proposition when $\nu_1=1$, using some implicit Presburger arithmetic on
the value group. We can partition $E$ into parts $E_1$ and $E_2$, with
\begin{eqnarray*}
 E_1 & = & \{x\in
 E| v(x_m)< v(\beta\prod_{i=2}^{m-1}x_i^{\nu_i})\},\\
 E_2 & = & \{x\in
 E| v(\beta\prod_{i=2}^{m-1}x_i^{\nu_i})\leq v(x_m)\}\\
 & = & \{x\in\prod_{i=1}^m R^{(\bar k)}| v(\beta\prod_{i=2}^{m-1}x_i^{\nu_i})\leq
  v(x_m)\leq v(\beta x_1\prod_{i=2}^{m-1}x_i^{\nu_i})\}.
\end{eqnarray*}
Since $ v(\beta\prod_{i=2}^{m-1}x_i^{\nu_i})\leq
 v(x_1\beta\prod_{i=2}^{m-1}x_i^{\nu_i})$ for $x\in E_1$, it follows that
 \[
 E_1=R^{(\bar k)}\times\{(x_2,\ldots,x_m)\in\prod_{i=2}^mR^{(\bar k)}| v(x_m)<
 v(\beta\prod_{i=2}^{m-1}x_i^{\nu_i})\},
 \]
and the restrictions $b_j\circ f|E_1$ satisfy condition
(\ref{condition}).\\
 As for $E_2$, let $D_{m-1}$ be the set
\[
D_{m-1}=\{(x_2,\ldots,x_m)\in\prod_{i=2}^mR^{(\bar k)}|v
 (\beta\prod_{i=2}^{m-1}x_i^{\nu_i})\leq v(x_m)\}.
\]
We may suppose that $\beta\in K^{(\bar k)}$, then the map
\[
R^{(\bar k)}\times D_{m-1}\to E_2:
x\mapsto(\frac{x_1x_m}{\beta\prod_{i=2}^{m-1}x_i^{\nu_i}},x_2,\ldots,x_m),
\]
can be checked by elementary Presburger arithmetic to be an isomorphism.
This proves the claim when $\nu_1=1$.
\par
Suppose now that $X$ is of the form described in (\ref{eq dim>1}) and
$\nu_1>1$. We prove we can reduce to the case $\nu_1=1$ by partitioning
and applying appropriate power maps. Choose $\tilde k>v(\nu_1)$ and put
$\tilde k'=\tilde k+v(\nu_1)$. We may suppose that $\tilde k\geq\bar k$,
so we have a finite partition $E=\bigcup_\alpha E_\alpha$, with
$\alpha=(\alpha_1,\ldots,\alpha_m)\in K^m$, $v(\alpha_1)=0$, $0\leq
v(\alpha_i)<\nu_1$ for $i=2,\ldots,m$ and
\[E_\alpha=\{x\in E|x_1\in \alpha_1R^{(\tilde k)},\ x_i\in \alpha_i
P_{\nu_1}^{(\tilde k')} \mbox{ for } i=2,\ldots,m\}.\] By corollary
\ref{corhensel} we have isomorphisms
\[f_\alpha:C_\alpha\to E_\alpha:x\mapsto(\alpha_1x_1,\alpha_2x_2^{\nu_1},\ldots,
\alpha_mx_m^{\nu_1}),\]
 with $C_\alpha=\{x\in\prod_{i=1}^m R^{(\tilde k)}|v(x_m)\leq v
(\beta' x_1\prod_{i=2}^{m-1}x_i^{\nu_i})\},$ where $\beta'\in K^\times$
depends on $\alpha$. This reduces the problem to the case described in
(\ref{eq dim>1}) with $\nu_1=1$ and thus the proposition is proved.
\end{proof}
\begin{remark} Proposition~\ref{resolution} can be strengthened, adding conditions on
the Jacobians of the isomorphisms, to become useful to $p$-adic
integration theory.\end{remark}
\begin{theorem}\label{dim}
Let $X$ be a semi-algebraic set, then either $X$ is finite or there exists
a semi-algebraic bijection $X\to K^k$ with $k\in\N_0$ the dimension of
$X$.
\end{theorem}
\begin{proof} We give a proof by induction on dim$(X)=m$. Let dim$(X)=1$. Use
Proposition \ref{resolution} to partition $X$ such that each part is
isomorphic to either $R^{(k)}$ or a point. By combining the isomorphisms
of Proposition~\ref{gring} and Corollary~\ref{corgring}, it follows that
$X\cong K$.
\par
 Now suppose dim$(X)=m>1$. Proposition 3 together with the case $m=1$
implies that we can finitely partition $X$ such that each part is
isomorphic to $K^l$, for some $l\in\{0,\ldots,m\}$, with
$K^0=\{0\}$. By proposition 2 at least one part must be isomorphic
to $K^m$. Suppose that $A$ and $B$ are disjoint parts, such that
$A\cong K^l$ and $B\cong K^m$, with $l\in\{0,\ldots,m\}$. It is
enough to prove that $A\cup B\cong K^m$. First suppose that $l=0$,
so $A$ is a singleton $\{a\}$. Since $m>1$ there exists an
injective semi-algebraic function $i:R\to A\cup B$ such that
$i(R\setminus\{0\})\subset B$ and $i(0)=a$. It follows that $A\cup
B\cong B\cong K^m$ since $R\cong R\setminus\{0\}$ (Proposition~1).
If $1\leq l$ we have $A\cup B \cong K\times(A'\cup B')$, for some
disjoint sets $A'\cong K^{l-1}$ and $B'\cong K^{m-1}$.  By
induction we find $A'\cup B'\cong K^{m-1}$ and thus $A\cup B\cong
K^m$. This proves Theorem 2.
\end{proof}
We obtain as a corollary of Theorem \ref{dim} the following
classification of the $p$-adic semi-algebraic sets.
\begin{cor}
Two infinite semi-algebraic sets are isomorphic if and only if they have
the same dimension.
\end{cor}
\begin{remark}
The field $\F_q((t))$ of Laurent series over the finite field is
often considered to be the characteristic $p$ counterpart of the
$p$-adic numbers. D. Haskell and the author \cite{CH} proved that
the Grothendieck ring of $\F_q((t))$ in the language of rings is
trivial, analogously as for $\Q_p$. In \cite{CH}, proof of Thm.~2,
it is also shown that there is a definable injection from the
plane $\F_q((t))^2$ into the line $\F_q((t))$, which makes it
plausible we cannot define an invariant for this field which
behaves like a (good) dimension. It is an open question to
classify the definable sets of $\F_q((t))$ up to definable
bijection. \end{remark}

 Raf Cluckers,\\
Katholieke Universiteit Leuven\\ Celestijnenlaan 200B\\ B-3001
Leuven,
Belgium;\\
 e-mail address: raf.cluckers@wis.kuleuven.ac.be
\end{document}